\documentclass[leqno,11pt,a4paper]{amsart}
\usepackage{pdfsync}
\usepackage{amsmath}
\usepackage{amsfonts}
\usepackage{amssymb}
\usepackage{tgbonum}
\usepackage{amsthm}
\usepackage{fullpage}
\usepackage{hyperref}
\usepackage{graphicx,multirow}
\usepackage{xcolor}
\usepackage{mwe}
\usepackage{tikz}
\usepackage{tikz-cd}
\hypersetup{
 colorlinks,
 linkcolor={blue!50!black},
 citecolor={blue!50!black},
 urlcolor={blue!80!black}
}
\usetikzlibrary{circuits.logic.US,circuits.logic.IEC,fit}

\usepackage{arydshln}
\usepackage{multirow}
\usepackage[all]{xy}
\usepackage{mleftright}
\usepackage{enumerate}
\usepackage{url}

\usepackage{pdflscape}
\usepackage{longtable}
\usepackage{booktabs}
\usepackage{colortbl}
\newcommand{\evnrow}{\rowcolor[gray]{0.95}}
\newcommand{\oddrow}{}

\usepackage{array,etoolbox}
\preto\tabular{\setcounter{magicrownumbers}{0}}
\newcounter{magicrownumbers}

\renewcommand{\P}{\mathbb P}

\newcommand{\Q}{\mathbb Q}
\newcommand{\M}{\mathbb M}
\newcommand{\Pf}{\mathbb{PF}}
\newcommand{\Z}{\mathbb Z}

\newcommand{\V}{\mathbb V}

\newcommand{\CI}{\mathbb{CI}}
\newcommand{\HS}{\mathbb{HS}}

\newcommand{\into}{\hookrightarrow}

\DeclareMathOperator{\Gr}{Gr}

\DeclareMathOperator{\Tom}{Tom}

\newcommand{\PxP}{\P^2\times\P^2}
\newcommand{\PI}{\P^1\times\P^1 \times\P^1}

\newtheorem{thm}{Theorem}

\theoremstyle{definition}
\newtheorem{dfn}[thm]{Definition}
\theoremstyle{remark}
\newtheorem{rmk}[thm]{Remark}
\theoremstyle{remark}

\numberwithin{equation}{section}
\numberwithin{thm}{section}

\newcommand{\QED}{\ifhmode\unskip\nobreak\fi\quad\ensuremath{\mathrm{QED}}}

\title{Projection Cascades of  models of log del Pezzo surfaces}
\author{Muhammad Imran Qureshi}
\date{}
\address{Muhammad Imran Qureshi\newline
Department of Mathematics, King Fahd University of Petroleum and Minerals (KFUPM), Dhahran
31261, Saudi Arabia \&\newline
Interdisciplinary Research Center for Intelligent Secure Systems, King Fahd University of Petroleum \& Minerals (KFUPM), Dhahran
31261, Saudi Arabia}
\email{imran.qureshi@kfupm.edu.sa}
\begin{document}

\maketitle   

\begin{abstract}
   We   introduce the notion of   type-I projection cascade for a biregular
 model (infinite series) of  log del Pezzo surfaces. We study the existence of type-I projection cascades for   known classes of models of  log del Pezzo surfaces, such that their images under their anti-canonical  embeddings in some weighted projective space, can be described as codimension 4 and codimension 3  varieties.  We obtain two cascades of length three and four cascades of length two, where each projection gives rise to a well-formed and quasismooth biregular model in lower codimension.    

\end{abstract}

%
\section{Introduction}
\subsection{Cascades  of log del Pezzo surfaces and graded rings} An algebraic surface $X$ with an ample anti-canonical class is called  the del Pezzo surface. They are  among one of the most studied mathematical structures in algebraic geometry.  There are a total of 10 deformation families of smooth del Pezzo surfaces: \(\P^1\times \P^1, \P^2,\) and blow ups of \(\P^2\) in \(d\) general points for \(1\le d \le 8:   \) forming a so called cascade of smooth  del Pezzo surfaces. A singular del Pezzo surface is
usually called a log del Pezzo surface or an orbifold del Pezzo surface if
the singularities are cyclic quotient singularities; we use both terms interchangeably. The construction and classification
of log del Pezzo surfaces have been the subject of great interest over the years from various perspectives, see \cite{Jhonson-Kollar,dp-exceptional,DP-zoo,dp-CH,Fano-CK,QMOC,dp-CP,QMATH}.

A log  del Pezzo surface with large
anti-canonical linear system  can be blown up successively to obtain a chain
of del Pezzo surfaces, called  a cascade of del Pezzo surfaces \cite{RS}.   
The blow up can be   considered  as a projection from a weighted projective space of a higher dimension to  a weighted projective space of a lower dimension. In favorable cases, this leads to  explicit descriptions of the images of
log del Pezzo surfaces under their anti-canonical embeddings  in some weighted projective space, as low codimension Gorenstein varieties. 

The cascade of log del Pezzo surfaces obtained by 
 blowing up   \(\P(1,1,3)\) in \(1\le d\le 8\) points was studied by Reid and Suzuki in \cite{RS}. They explicitly described the Gorenstein construction
of these varieties  in low codimension cases. In particular, they showed  that del the Pezzo surfaces  for $d$ equal to \(6,7,\)  and \(8\) have descriptions  in terms of codimension 3 variety in \(\P^5(w_i)\) defined by the maximal Pfaffians of \(5 \times 5
\) skew-symmetric matrix of general forms , codimension 2 complete intersection variety in \(\P^4(w_i)\) and a hypersurface in \(\P^3(w_i)\) respectively.  Moreover, they are connected by  type-I  Gorenstein projection \cite{RP}
 from  a smooth coordinate point of  the ambient weighted projective space \(\P(w_i),\) followed by a   deformation of the image of projection to a generic member of some
family, forming what we call   a type-I projection cascade for log del Pezzo
surfaces, see Definition \ref{def:cascade}. Also,
the case \(d=5\) can  be described as a codimension 4 variety that can
be obtained as an unprojection of the Pfaffian del Pezzo surface that corresponds
to   \(d=6\).  

\subsection{Aim and results}  

The aim of this paper is to study the type-I projection cascades for biregular models
of log del Pezzo surfaces,  introduced in \cite{QMOC}, obtained via algorithmic
developed in  \cite{QJSC}. For this purpose we
introduce the  notion of type-I projection cascade for biregular
models of log del Pezzo surfaces and  study if our biregular models from \cite{QMOC,QP13} may form a type-I projection cascade. A  biregular model of log del Pezzo surfaces is an infinite series of such
surfaces   parameterized by the positive integers, formally defined
as follows.
\begin{dfn}\cite{QMOC}
\label{def:models}
A \emph{biregular model of log del Pezzo surfaces} \(\M\) is  an infinite series of families of   del Pezzo surfaces satisfying the following conditions.  
\begin{enumerate}
\item  A family of  wellformed and quasismooth  log del Pezzo surfaces exists for each  parameter $r(n)$ for all $n \in \Z_{\ge 0}$, where \(r \) is a linear function of
\(n\). 
\item The
ambient weighted projective space \(\P(w_i) \) that contains the image of  log del Pezzo surface under their (sub) anti-canonical embedding must
 contain at least
one weight to be  $r$, and  rest of the weights  $w_i$  and the volume  $(-K_X)^2$ are  functions of $r$.
\end{enumerate}
\end{dfn}
A simple example of a biregular model is an index \(r+2\) del Pezzo surface
\(\P(1,1,r)\) that contains an orbifold point of type \(\frac1r(1,1)\),
where \(r=n \in \Z_+\).  
\noindent
   
\begin{dfn}
\label{def:cascade}A \emph{ type-I projection cascade of biregular models of log  del
Pezzo surfaces} is a collection   \(\{\M_1,\ldots,\M_k\}\) of biregular models of log del Pezzo surface such that for any  surface \(X_1^{r(n)}\) in \(\M_1\) there exist a collection of surfaces \(\{X_i^{r(n)}: X_i^{r(n)} \in \M_i, 2\le i\le k\}\) such that   
\begin{enumerate}
\item  \(X_1^{r(n)}\) has image under their  anti-canonical  embedding  described as a subvariety of \(\P(1^l, a_1, a_2, \ldots, a_s)\) where \(l \ge k\) and \(a_i\ge 2 \), and each subsequent
del Pezzo surface \(X_i^{r(n)}\)  under their  anti-canonical embedding is a subvariety of \(\P(1^{l-i+1}, a_1, a_2, \ldots, a_s)\), i.e. \(X_{i+1}^{r(n)}\) has one codimension lower than \(X_{i}^{r(n)}. \) 
\item  There exist type-I projection from a smooth coordinate point \(p_{0}\) of weight one  from \(X_{i}^{r(n)}\)  to \(\overline X_{i+1}^{r(n)}\supset D \), followed
by a deformation of \(\overline X_{i+1}^{r(n)}\) to a wellformed and quasismooth log del Pezzo surface  \(X_{i+1 }^{r(n)}\in \M_{i+1}.\) The divisor \(D\) is the image of the smooth coordinate point \(p_{0}\) of the variable  \(x_0\) under the projection map \(\pi_i\).  
\[
\begin{tikzcd}[row sep=2em]
p_{0}\in X_i^{r(n)}\arrow[dashed,rightarrow]{rd}{\pi_i}\arrow[rightarrow]{rr}{\phi_i}  & &  X_{i+1}^{r(n)}\\
&D\subset \overline X_{i+1}^{r(n)}\arrow[hookrightarrow]{ru}{\psi_i} &
\end{tikzcd}
\]
\end{enumerate}
One can alternatively rephrase this definition to introduce an unprojection cascade,
in particular the condition   $(ii)$ can be restated as: There exists a degeneration
of \(X_{i+1}^{r(n)}\rightsquigarrow \overline X_{i+1}^{r(n)}\supset D\), followed by  unprojection to  \(X_i\) which contracts the divisor \(D\) to the coordinate point  $p_{_0}$ of the variable \(x_0\).  
\end{dfn}
 
 It follows from the definition that the existence of type-I
projection cascade is possible for those models where the ambient weighted
projective space \(\P(w_i)\) contains at least one variable of  weight 1, so we list only such models
from \cite{QMOC} in Appendix \ref{list}. In total, there are 9 such models but only 6 of them give type-I projection cascades for the models of log del Pezzo surfaces.  Projections from the rest of the three models do not deform to a quasismooth
model of log del Pezzo surfaces. Moreover, no models from \cite{QP13} fit into a type-I projection
cascades of models of log del Pezzo surfaces.   

In \cite{QMOC,QP13}, we constructed several models of rigid log del Pezzo surfaces whose equations were given in terms of equations  of  either
 the Pl\"ucker embedding of Grassmannians $\Gr(2,5)  $   or of the equations of the Segre     embedding of \(\PxP\) in \(\P^8\) or of  the Segre     embedding of \(\PI\) in \(\P^7\), giving models as codimension 3, 4 and 4 Gorenstein varieties in some weighted projective
space,  respectively.
We get the following result, which requires constructions  of some codimension 2 complete intersection models and hypersurface models to complete the proof.

 \begin{thm}
\label{thm:Cascade}
Let \(\P_{ij}\) denote a codimension 4 \(\mathbb P^2 \times \mathbb P^2\)
model, \(\Pf_{ij}\) denote a codimension  3 Pfaffian model, \(\CI_{ij}\) be
a complete intersection model and \(\HS_{ij}\) be  a hypersurface model  of log del Pezzo
surfaces then there exist six type 1 projection cascades of    models of log del Pezzo surfaces  given as follows. 
\begin{enumerate}[(I)]

\item 
\begin{tikzcd}[row sep=2em]
\P_{11}\arrow[rightarrow]{rr}  & &   \Pf_{11} \arrow[rightarrow]{rr}  & &   \CI_{11} 
\end{tikzcd}

 \item 
 \begin{tikzcd}[row sep=2em] 
 \Pf_{12} \arrow[rightarrow]{rr}  & &   \CI_{12}\arrow[rightarrow]{rr}  & &    \HS_{12}
\end{tikzcd}

 \item \begin{tikzcd}[row sep=2em] 
 \Pf_{13} \arrow[rightarrow]{rr}  & &   \CI_{13} 
\end{tikzcd}
 \item \begin{tikzcd}[row sep=2em] 
 \Pf_{21} \arrow[rightarrow]{rr}  & &   \CI_{21} 
\end{tikzcd}
 \item \begin{tikzcd}[row sep=2em] 
 \Pf_{22} \arrow[rightarrow]{rr}  & &   \CI_{22} 
\end{tikzcd}

\item \begin{tikzcd}[row sep=2em] 
\P_{12} \arrow[rightarrow]{rr}  & &   \Pf_{14} 
\end{tikzcd}
\end{enumerate}
where the  descriptions of complete intersection and hypersurface models are listed in the following table and the rest of the models are given in Appendix
\ref{list}.  
\end{thm}
 \renewcommand*{\arraystretch}{1.2}
\begin{longtable}{>{\hspace{0.5em}}llccccr<{\hspace{0.5em}}}
\caption{Complete intersection and Hypersurface  Models,       $q=r-1, s=r+1,t=r+2,u=r+3,v=2r+1,y=2r-2,z=2r-1,m=3r-2$. } \label{tab-CI}\\
\toprule
\multicolumn{1}{c}{Model}&\multicolumn{1}{c}{WPS  }&\multicolumn{1}{c}{Basket $\mathcal{B}$}&\multicolumn{1}{c}{$-K^2$}&\multicolumn{1}{c}{$h^0(-K)$}\\
\cmidrule(lr){1-1}\cmidrule(lr){2-2}\cmidrule(lr){3-3}\cmidrule(lr){4-4}\cmidrule(lr){5-5}
\endfirsthead
\multicolumn{7}{l}{\vspace{-0.25em}\scriptsize\emph{\tablename\ \thetable{} continued from previous page}}\\
\midrule
\endhead
\multicolumn{7}{r}{\scriptsize\emph{Continued on next page}}\\
\endfoot
\bottomrule
\endlastfoot
\evnrow $\CI_{11}$ & 
   $ \begin{array}{@{}l@{}} Y_{2r,2r}\\\quad\subset\P(1^2,r^2,z);\\r=n+1, n\ge 1\end{array}$ & $\frac{1}{z}(1,1)$&$\frac{4}{2 r-1}
$&$  2 $\\ 

 \oddrow  $\CI_{12}$ &  $ \begin{array}{@{}l@{}} Y_{2r,2r+1}\\\quad\subset\P(1,2 ,r^2,z);\\r=2n+1, n\ge 1\end{array}$ & $2\times\frac{1}{r}(2,q),\frac{1}{z}(1,1)$&$\frac{2 r+1}{r(2 r-1)}
$&$  1 $\\  
\evnrow $\HS_{12} $& $ \begin{array}{@{}l@{}} Z_{4r}\\\quad\subset\P(2 ,r^2,z);\\r=2n+1, n\ge 1\end{array}$ & $4\times\frac{1}{r}(2,q),\frac{1}{z}(1,1)$&$\frac{2}{r (2 r-1)}
$&$  0 $\\  

\oddrow $\CI_{13} $&$ \begin{array}{@{}l@{}} Y_{3r,4r-2}\\\quad\subset\P(2 ,r^2,z,m);\\r=2n+1, n\ge 1\end{array}$ & $3\times\frac{1}{r}(2,q),\frac{1}{m}(r,z)$&$\frac{3}{r (3r-2)}
$&$  0 $\\

\evnrow $\CI_{21} $&$ \begin{array}{@{}l@{}} Y_{r+2,r+3}\\\quad\subset\P(1,2,3,r,s);\\r=3n, n\ge 1\end{array}$ & $
    \begin{matrix}
\frac{1}{2}(1,1),\frac{1}{3}(1,1),\\\frac{1}{r}(1,1),\frac{1}{s}(3,r)
\end{matrix}$&$\frac{(r+2) (r+3)}{6 r (r+1)}
               $&$  2 $\\  
\oddrow $\CI_{22} $&$ \begin{array}{@{}l@{}} Y_{r+2,r+3}\\\quad\subset\P(1,2,3,r,s);\\r=3n+1, n\ge 1\end{array}$ & $
    \begin{matrix}
\frac{1}{2}(1,1),\frac{1}{r}(1,1),\\\frac{1}{s}(3,r)
\end{matrix}$& Same as \(\CI_{21}\)&$  2 $\\   
  
\end{longtable}  
\begin{rmk} The cascade \((I)\) was given in \cite{dp-CP} but the description  in terms of equations of surfaces in the model \(\CI_{11}\)  did not appear.  
\end{rmk}

 
\section{Reid--Suzuki cascade and graded rings}
 As
a warm up example, we  recall  the Ried--Suzuki cascade of del Pezzo surface obtained from \(\P(1,1,3) \) and its relation with type-I projection and  graded rings from \cite{RS}.  The blow up of \(k\) general points \(P_i\) of \(\P(1,1,3)\)  with \(k\le 8\) forms a cascade of log del Pezzo surface. If \(S^{(k)}\to \P(1,1,3)\) is the  blow up  in \(k \) general points, then for \(5\le k\le 8\)   explicit graded ring constructions  for these surfaces  have been described as follows.   

\begin{itemize}
\item 
\(k=8\):  has a description as a hypersurface \(X_{10}\into \P(1,2,3,5).\) 

\item 
\(k=7\): can be described as a complete intersection of two weighted homogeneous quartics \(X_{4,4}\into \P(1,1,2,2,3)\).
\item 
\(k=6\):  a maximal Pfaffians of   \(5\times 5\) skew symmetric matrix of  weighted homogeneous forms   in the graded ring of    \( \P(1^{3},2^{2},3)=\P(x_0,x_1,x_2,y_0,y_1,z)\). The Pfaffian matrix is given below where \(Q_{i}\) represents generic quadrics.  We omit the zeros at the diagonal and lower triangular part of the matrix   from the presentation. 
\[\begin{pmatrix}x_0&x_1&Q_1&Q_2\\&x_2&Q_3&Q_4\\&&y_0&y_1\\&&&z \end{pmatrix}\textrm{ of degrees}\begin{pmatrix}1&1&2&2\\&1&2&2\\&&2&2\\&&&3 \end{pmatrix}. \]
\item 
\(k=5\):  a codimension 4 subvariety in \(\P(1^4,2^2,3) \) obtained as a Kustin--Miller unprojection of \(k=6\) case  by inserting a divisor \(\P^1\) in \(S^{(6)}\).
The graded ring can be explicitly calculated and it has a \(9 \times 16 \) resolution; though they are not known to be in a specific equation format.   \end{itemize}

 We can think of the surface \(S^{(k)}\) as a type-I unprojection \cite{RP}  of the degeneration \(\overline S^{(k-1)}\) of \(S^{(k-1)}\), containing the curve \(C=\P(1,1)\) for cases where \(k \) is equal to  \(5,6\) and \(7\). Equivalently, we can observe that the cases for \(k=5,6 \) and \(7\) are linked via type-I projection and they fit into a projection cascade.   For example, the complete intersection \(X_{4,4}\) can be obtained as a type-I projection from the point \(x_0\) of the codimension 3 Pfaffian variety, followed by flat deformation to  a generic complete intersection family. Indeed the projection itself will give us \({\overline X}_{4,4}\) containing the unprojection divisor \(D=\P(1,1)\) and it is a flat degeneration of  general  complete intersection \(X_{4,4}\into\P(1^2,2^2,3)\).


\section{Cascades of projections of log del Pezzo surfaces}
%

\subsection{Proof of Theorem \ref{thm:Cascade} }
 In this section, we prove the main theorem by describing the constructions
in all cases. Throughout the proof \(S, X, Y\) and \(Z\) denote  del Pezzo surfaces of  codimension
$4,3,2$ and $1$ respectively, in the given biregular models. The symbols
\(F_{j}, G_j, H_j,\) and \(I_j\) denote generic weighted homogeneous forms of degree \(j\) and \(x_0,y_0, y_{1, }z_{0, }a,b,c,d,e\) and \(f\) denote the variables of the ambient
weighted projective space.    
\subsubsection{Case I} In the model \(\P_{11}\) the image of each del Pezzo surface
 \(S\) under their anti-canonical embedding can be described as a codimension 4 subvariety
 of \[\P(1^4,r^2,2r-1)=\P(x_0,y_{0,}a,b,c,d,e). \]
Following Appendix \ref{list}, a set of equations defining image of   \(X\) are  given  by the \(2\times 2\) minors of the  matrix, 

\[\begin{pmatrix}
x_0 &y_0 &c\\
a&b&d\\
F_{r}&G_{r}&e
\end{pmatrix},\]
 where \(F_{r}\) and \(G_{r}\) are weighted homogeneous forms of degree \(r\).
 
 The implicit functions near the point \(p_{x_0}\) are \(b,d,e\) and \(G_{2r}\) and the image of Gorenstein projection \(\widetilde X\) from \(p_{x_0}\) is given by the maximal Pfaffians of the following skew symmetric matrix,
where we omit the diagonal and lower triangular parts in the presentation. 
\[N= \begin{pmatrix}
y_{0}&c&a&  F_{r}\\
&   0& b       &d\\
&& G_{r}&e\\&&& 0
\end{pmatrix}. \] The surface \(\widetilde X\subset \P(1^{3},r^2,2r-1)\) contains the image of \(p_{x_0}\) to be the divisor \[D=\V(b, d,e,G_{r}).\] The matrix \(N\) is in the \(\Tom_1\) format \cite{BKR}, that is, all entries of \(N\), except the first row and column, are in the ideal of the divisor \(D\)  for  \(I_D:=<b,d,e,G_{2r}>\).   The zero polynomials at position \((2,3)\) and  $(4,5)$ are of degrees \(r\) in the matrix \(N\). 

The \(\Q\)-smoothing of \(\widetilde X\) to \(X\) is obtained by deforming the entries of \(N\) to be generic polynomials in the variables of the ambient weighted projective
space, which gives us the  model \(\Pf_{11}\).  Each  \(X\) in the model
\(\Pf_{11}\) under their anti-canonical  embedding can be described as a
subvariety of  \[\P(1^3,r^2,z)=\P(y_{0,}a,b,c,d,e).\] 
  Its equations are given by the     
maximal Pfaffians of the \(5\times 5 \) skew symmetric matrix  
\[M= \begin{pmatrix}
y_0&c&a&  F_r\\
&   H_{r}& b       &d\\
&& G_r&e\\&&& I_{r}
\end{pmatrix},\]
where \(I_r\) and \(H_r\) are generic forms and we omit the diagonal and lower triangular parts. 

Now we perform a Gorenstein projection from the point \(y_{0}\).
It is clear from the equations that \(G_{r}, e\) and \(I_{r}\) are implicit functions around \(p_{y_0}\). Therefore the image of \(X\) under the Gorenstein projection from \(p_{y_0}\) is a complete intersection \(\widetilde Y\): 
\[\widetilde Y_{2r,2r}:\V\left(\left[\begin{matrix}c&  a&F_r\\H_r& b       &d\end{matrix}\right]\left[\begin{matrix}I_{r} \\e\\G_r\end{matrix}\right] \right)\subset \P(1^{2},r^2,2r-1). \]  
The divisor \(D:= \V(I_{r},G_{r+1},e,)\subset \widetilde Y\) is the image of the point \(P_{y_0}\) under the Gorenstein projection. Then the \(\Q\)-smoothing
 \(\widetilde Y\to Y\) gives a generic complete intersection of degrees \(2r\):  \[Y_{2r,2r}\subset \P(1^2,r^2,2r-1)=\P(a,b,c,d,e), \] 
which is a wellformed and quasismooth model of log del Pezzo surfaces described as \(\CI_{11}\) in Table \ref{tab-CI}.
 
 There are two further variables of weight one, however \(\Q\)-smoothing
 after Gorenstein projection does not lead to a quasismooth hypersurface
 model  \(Z_{4r-1}\subset \P(1,r^2,2r-1).  \)  
\subsubsection{Case II} The biregular model \(\Pf_{12}\) consists of log del Pezzo surfaces that can be embedded in \[\P(1^2,2,r^2,2r-1)=\P( y_0,z_0,a,b,c,d). \] Following Appendix \ref{list}, the equations of any quasismooth surface \(X\) in this model are given by  \( 4\times 4 \) Pfaffians of the \(5\times 5 \) skew symmetric matrix   
\[M=
 \begin{pmatrix}
y_0&z_0&F_{r-1}&  F_r\\
&   a& b       &F_{r+1}\\
&& c&G_{r+1}\\&&& d
\end{pmatrix},\] where we omit the diagonal and lower triangular parts.   
It is clear from the equations that \(c, d\) and \(G_{r+1}\) are implicit functions around \(p_{y_0}\). Therefore the image of \(X\) under the Gorenstein projection from \(p_{y_0}\) is a complete intersection \(\widetilde Y\): 
\[\widetilde Y_{2r,2r+1}:\V\left(\left[\begin{matrix}z_0&F_{r-1}&  F_r \\a& b       &F_{r+1}\end{matrix}\right]\left[\begin{matrix}d \\G_{r+1}\\c\end{matrix}\right] \right)\subset \P(1,2,r^2,2r-1). \]  
The divisor \(D:= \V(G_{r+1},c,d)\subset \widetilde Y\) is the image of the point \(p_{y_0}\) under the Gorenstein projection. The \(\Q\) smoothing  \(\widetilde Y \to Y\)  is generic complete intersection of degrees \(2r\) and \(2r+1\):  \[Y_{2r,2r+1}\subset \P(1,2,r^2,2r-1)=\P(z_0,a,b,c,d), \] 
which is a wellformed and quasismooth model of log del Pezzo surfaces described as \(\CI_{12}\) in Table \ref{tab-CI}. 

Now we perform a Gorenstein projection from the point \(z_{0}\) which indeed amounts to the elimination of variable \(z_0\).  We can write the equations of \(Y \) as 
\[z_0d= H_{2r}(a,b,c,d),\; z_0G_{2r}=H_{2r+1}(a,b,c,d).\] 
The image of \(Y\) under the projection from \(z_0\) is a hypersurface \[\widetilde Z_{4r}\subset \P(2,r^2,2r-1),\]
containing the divisor \(E=\V(d,G_{2r})\). The \(\Q\)-smoothing \( \widetilde
Z \to Z \) is the generic hypersurface of degree \(4r\) in \(\P(2,r^2,2r-1)\), giving rise to a  wellformed and quasismooth biregular model \(\HS_{12}\) given in the Table \ref{tab-CI}. 
\subsubsection{Case III} The biregular model \(\Pf_{13}\) consists of log del Pezzo surfaces that can be embedded in \[\P(1^2,2,r^2,2r-1,3r-2)=\P( y_0,y_{1},a,b,c,d,e). \] Following Appendix \ref{list}, the equations of any quasismooth surface \(X\) in this models can be described by  \( 4\times 4 \) Pfaffians of the \(5\times 5 \) skew symmetric matrix   
\[M=
 \begin{pmatrix}
y_0&a&b&  F_{r+1}\\
&   c& F_{2r-2}       &d\\
&& F_{2r-1}&F_{2r}\\&&& e
\end{pmatrix},\] where \(F_i\)s denote generic forms of degree \(i\).   
It is clear from the equations that \(e, F_{2r}\) and \(F_{2r-1}\) are implicit functions around \(p_{y_0}\). Therefore the image of \(X\) under the Gorenstein projection from \(p_{y_0}\) is a complete intersection \(\widetilde Y\): 
\[\widetilde Y_{4r-2,3r}:\V\left(\left[
\begin{matrix}a&b&  F_{r+1} \\
c& F_{2r-2}       &d\end{matrix}\right]\left[\begin{matrix}e \\F_{2r}\\F_{2r-1}\end{matrix}\right] \right)\subset \P(1,2,r^2,2r-1,3r-2). \]  
The divisor \(D:= \V(e,F_{2r},F_{2r-1})\subset \widetilde Y\) is the image of the point \(p_{y_0}\) under the Gorenstein projection. The \(\Q\)-smoothing
of \(\widetilde Y \to Y\)  is generic complete intersection of degrees \(2r\) and \(2r+1\):  \[Y_{4r-2,3r}\subset \P(1,2,r^2,2r-1)=\P(y_0,a,b,c,d), \] 
which is a wellformed and quasismooth model of log del Pezzo surfaces, given as \(\CI_{13}\) in Table \ref{tab-CI}. 
\subsubsection{Case IV}The biregular model \(\Pf_{21}\) consists of log del Pezzo surfaces that can be described under their anti-canonical embeddings
as subvarieties of  \[\P(1^2,2,3,r,r+1)=\P( y_0,a,b,c,d,e). \] Following Appendix \ref{list}, the equations of any quasismooth surface \(X\) in this model can be described by  \( 4\times 4 \) Pfaffians of the \(5\times 5 \) skew symmetric matrix   
\[M=
 \begin{pmatrix}
y_0&a&b&  F_{r-1}\\
&   F_{2}& c&d\\
&& F_{3}&F_{r}\\&&& e
\end{pmatrix},\] where \(F_i\)s denote generic forms of degree \(i\).   
It is clear from the equations that \(e, F_{r}\) and \(F_{3}\) are implicit functions around \(p_{y_0}\). Therefore the image of \(X\) under the Gorenstein projection from \(p_{y_0}\) is a complete intersection \(\widetilde Y\): 
\[\widetilde Y_{r+2,r+3}:\V\left(\left[
\begin{matrix}a&b&  F_{r-1} \\
F_{2}& c       &d\end{matrix}\right]\left[\begin{matrix}e \\F_{r}\\F_{3}\end{matrix}\right] \right)\subset \P(1,2,3,r,r+1). \]  
The divisor \(D:= \V(e,F_{r},F_{3})\subset \widetilde Y\) is the image of the point \(p_{y_0}\) under the Gorenstein projection. The \(\Q\)-smoothing
of \(\widetilde Y \to Y\)  is generic complete intersection of degrees \(r+2\) and \(r+3\):  \[Y_{r+2,r+3}\subset \P(1,2,3,r,r+1)=\P(a,b,c,d,e), \] 
which is a wellformed and quasismooth model of log del Pezzo surfaces described as \(\CI_{21}\) in Table \ref{tab-CI}. 
\subsubsection{Case V} This case is exactly similar to case $IV$, except
that the parameter in case $IV$ is \(r=3n\) and in this case it is \(r=3n+1\).

\subsubsection{Case VI} This is a one step cascade which starts at a codimension 4 \(\PxP\) model \(\P_{12}\) and admits a projection to \(\Pf_{14}\). In model \(\P_{12}\) each surface \(S\) under their anti-canonical  embedding can be described as a subvariety of \[\P(1,2,r^2,r+1,r+2,2r-1)=\P(x_0,a,b,c,d,e,f). \]
The equations of \(X\) are given by the \(2\times 2\) minors of following matrix, where \(F_{r+1}\) and \(F_{2r}\) are weighted homogeneous forms of degrees \(r+1\) and \(2r\), 
\[\begin{pmatrix}x_0 &a &b\\c&d&f\\F_{r+1}&e&F_{2r}\end{pmatrix}.\] The implicit functions near the point \(p_{x_0}\) are \(d,e,f\) and \(F_{2r}\) and the image of Gorenstein projection \(\widetilde X\) from \(P_{x_0}\) is given by the maximal Pfaffians of the following skew symmetric matrix. 
\[N= \begin{pmatrix}
a&b&c&  F_{r+1}\\
&   0& d       &f\\
&& e&F_{2r}\\&&& 0
\end{pmatrix}. \] The surface \(\widetilde X\subset \P(2,r^2,r+1,r+2,2r-1)\) contains the image of \(p_{x_0}\) to be the divisor \[D=\V(d,e,f,F_{2r})\cong\V(F_{2r})\subset \P(2,r^{2}).\] The matrix \(N\) is in the \(\Tom_1\) format \cite{BKR}, that is, all entries of \(N\), except the first row and column, are in the ideal of the divisor \(D\)  for  \(I_D:=<d,e,f,F_{2r}>\).  The divisor contains two \(\frac 1r\) type of points which contracts to the point \(p_{x_0}\) under the Gorenstein unprojection of \(\widetilde X\) to \(\PxP\) model \(\P_{12}\). The \(\Q\)-smoothing of \(\widetilde X\) to \(X\) is obtained by deforming the entries of \(N\) to be generic polynomials in the ambient variables, which gives us the  model \(\Pf_{14}\) of Appendix \ref{list}. The zero polynomials at position \((2,3)\) has degree \(r+1\)  and $(4,5)$ has degree \(2r\) in the matrix \(N\). There is no further weight one variable so we do not perform any further projections to a lower codimension.

\subsection*{Acknowledgment}

 I am grateful to Stephen Coughlan for helpful discussions on the subject
and to an anonymous referee for his/her comments that helped improve the exposition significantly. 
\bibliographystyle{amsalpha}
\bibliography{References}
\appendix\section{Tables of models of log del Pezzo surfaces}
\label{list}
We list the tables of models of log del Pezzo surfaces from \cite{QMOC} in this section. These tables contains those 9 biregular models of del Pezzo
surface that can be  embedded in  weights  weighted projective space  and
having at least one weight of degree one, i.e. those models where a cascade of projections possibly exists. 

 \label{thrm-Pf}  
 
 \renewcommand*{\arraystretch}{1.2}
\begin{longtable}{>{\hspace{0.5em}}llccccr<{\hspace{0.5em}}}
\caption{Biregular log del Pezzo Models in  Pfaffian format where  $q=r-1, s=r+1,t=r+2,u=r+3,v=2r+1,y=2r-2,z=2r-1,m=3r-2$. } \label{tab-Pf}\\
\toprule
\multicolumn{1}{c}{Model}&\multicolumn{1}{c}{WPS \& Param }&\multicolumn{1}{c}{Basket $\mathcal{B}$}&\multicolumn{1}{c}{$-K^2$}&\multicolumn{1}{c}{$h^0(-K)$}&\multicolumn{1}{c}{Weight Matrix}\\
\cmidrule(lr){1-1}\cmidrule(lr){2-2}\cmidrule(lr){3-3}\cmidrule(lr){4-5}\cmidrule(lr){6-6}
\endfirsthead
\multicolumn{7}{l}{\vspace{-0.25em}\scriptsize\emph{\tablename\ \thetable{} continued from previous page}}\\
\midrule
\endhead
\multicolumn{7}{r}{\scriptsize\emph{Continued on next page}}\\
\endfoot
\bottomrule
\endlastfoot
\evnrow $\Pf_{11}$ & 
  $\begin{matrix}
  \P(1^3, r^2, z)\\  \; r =n+1 \\ 
  w=\frac 12(1,1,1,z,z)
  \end{matrix}$&
    
$\dfrac{1}{z}(1,1)$
&  $\dfrac{2r+3}{2r-1} $             & 3&

$ \footnotesize{\begin{matrix} 1&1&r&r\\ &1&r&r\\ &&r&r\\ &&& z \end{matrix}}$  \\

\oddrow $\Pf_{12}$ & $\begin{matrix}
  \P(1^2, 2,r^2,z)\\  r=2n+1,\; n\ge 1 \\
  w=(0,1,1,q,r)  
  
  \end{matrix}$ & $\dfrac{1}{r}(2,q),\dfrac{1}{z}(1,1)$&$\dfrac{2 r^2+5 r+1}{4 r^2-2 r}$&$  2 $&$ \footnotesize{\begin{matrix} 1&1&q&r\\ &2&r&s\\ &&r&s\\ &&& z \end{matrix}} $\\  

\evnrow $\Pf_{13} $&$ 
  \begin{matrix}
  \P(1, 2,r^{2}, z,m),\\  r=2n+1,\; n \ge 1\\
 w= \frac 12(4-s,s-2,\\s,3s-6,3s-4)
  \end{matrix}$&$ 
    \begin{matrix}
 \dfrac{1}{r}(2,q), \dfrac{1}{m}(r,z)
\end{matrix}$&$\dfrac{3 r+1}{r (3 r-2)}$&$1
               $&$
 \footnotesize\begin{matrix} 1&2&r&s\\ &r&y&z\\ &&z&2r\\ &&& m \end{matrix} $ \\
\oddrow $\Pf_{14} $&$ 
  \begin{matrix}
  \P( 2,r^{2},s,t, z)\\  r=2n+1,\; n \ge 1\\
  w=\frac 12(1,3,z,z,2r+1)
  \end{matrix}$&$     \begin{matrix}
  \dfrac{1}{t}(2,s),\dfrac{1}{z}(1,1)\\3\times \dfrac{1}{r}(2,q),
\end{matrix}$&$\dfrac{4 r+3}{r (r+2) (2 r-1)}
               $&$0$&$

 \footnotesize\begin{matrix} 2&r&r&s\\ &s&s&t\\ &&z&2r\\ &&& 2r \end{matrix}$  \\

 \evnrow $\Pf_{21} $&$ 
  \begin{matrix}
  \P(1^2, 2,3,r, s);\\  r=3n,\; n \ge 2 \\
  w=(0,1,1,2,q)
  \end{matrix}$&$
    \begin{matrix}
\dfrac{1}{3}(1,1),\dfrac{1}{r}(1,1),\\\dfrac{1}{s}(3,r)
\end{matrix}$&$\dfrac{4 \left(2 r^2+4 r+3\right)}{3 r (r+1)}               $&$  4$&$

 \footnotesize{\begin{matrix} 1&1&2&q\\ &2&3&r\\ &&3&r\\ &&& s \end{matrix}}$ \\

\oddrow $\Pf_{22} $&$ \begin{matrix}
  \text{ rest same as } \Pf_{21}\\  r=3n+1,\; n \ge 2
  \end{matrix} $&$  \dfrac{1}{r}(1,1),\dfrac{1}{s}(3,r)$&$\text{ Same as } \Pf_{21} $&$ 4$&$ \text{Same as}  \Pf_{21}$\\  

\evnrow $\Pf_{23} $&$ 
  \begin{matrix}
  \P(1, 3,r, s,t,u)\\  r=3n+2,\; n \ge 2\\
  w=(0,1,2,r,s)
  \end{matrix}$&$
    \begin{matrix}
 \dfrac{1}{3}(1,1),\\ \dfrac{1}{r}(1,1),\dfrac{1}{u}(3,t)
\end{matrix}$&$\dfrac{8 r+36}{3 r^2+9 r}
               $&$1$&$
 \footnotesize\begin{matrix} 1&2&r&s\\ &3&s&t\\ &&t&u\\ &&& v \end{matrix} $ \\

\end{longtable}
  
 \renewcommand*{\arraystretch}{1.2}
\begin{longtable}{>{\hspace{0.5em}}llccccr<{\hspace{0.5em}}}
\caption{biregular log del Pezzo Models in  $\PxP$ format where $q=r-1, s=r+1,t=r+2, z=2r-1$. } \label{tab-P2xP2}\\
\toprule
\multicolumn{1}{c}{Model}&\multicolumn{1}{c}{WPS \& Param\ }&\multicolumn{1}{c}{$\mathcal{B}$}&\multicolumn{1}{c}{$-K^2_X$}&\multicolumn{1}{c}{$h^0(-K)$}&\multicolumn{1}{c}{Weight Matrix}\\
\cmidrule(lr){1-1}\cmidrule(lr){2-2}\cmidrule(lr){3-3}\cmidrule(lr){4-5}\cmidrule(lr){6-6}
\endfirsthead
\multicolumn{7}{l}{\vspace{-0.25em}\scriptsize\emph{\tablename\ \thetable{} continued from previous page}}\\
\midrule
\endhead
\multicolumn{7}{r}{\scriptsize\emph{Continued on next page}}\\
\endfoot
\bottomrule
\endlastfoot

\evnrow $\P_{11} $&$ 
  \begin{matrix}
  \P(1^4, r^2, z)\\  \; r=n+1\\w=(0,0,q; 1,1,r)
  \end{matrix}$&$

\dfrac{1}{z}(1,1)
$&$ \dfrac{4 r+2}{1-2 r}              $&$ 4$&$
 \footnotesize{\begin{matrix} 1&1&r\\ 1&1&r\\ r&r& z \end{matrix}} $ \\

\oddrow $\P_{12} $&$ 
  \begin{matrix}
  \P(1,2, r^2,s,t, z)\\  \; r=2n+1\; n\ge 1\\
  w=(0,1,q; 1,r,s)
  \end{matrix}$&$
\dfrac{1}{r}(2,q),\dfrac{1}{t}(2,s), \dfrac{1}{z}(1,1)
$&$\dfrac{2 r^2+7 r+1}{r (r+2) (2 r-1)}
 $&$ 1$&$  
 \footnotesize{\begin{matrix} 1&2&r\\ r&s&z\\ s&t& 2r \end{matrix}} $ \\

\evnrow $\P_{13} $&$ 
  \begin{matrix}
  \P( 1,2^2,3,r^{2}, z)\\  r=2n+1,\; n \ge 1\\
  w=(0,1,q;1,2,r)
  \end{matrix}$&$
    \begin{matrix}
 \dfrac{1}{3}(1,1), \dfrac{1}{z}(1,1)\\
 2\times \dfrac{1}{r}(2,q)
 
\end{matrix}$&$  \dfrac{(2 r+3)}{3 r (2 r-1)}
              $&$1 $&$

 \footnotesize\begin{matrix} 1&2&r\\ 2&3&s\\ r&s&z\end{matrix} $ \\

\end{longtable}
 

\end{document}